\documentclass[letterpaper, 10 pt, conference]{ieeeconf}
\usepackage{amssymb}
\usepackage{color}
\usepackage{bbm}

\usepackage{graphicx}
\usepackage{mathtools}
\usepackage{multirow}
\usepackage{lineno,hyperref}
\usepackage{url}
\usepackage{cleveref}
\usepackage{subcaption}
\usepackage{bm}
\usepackage{tikz}
\usepackage{algorithm}
\usepackage{algorithmicx}
\usepackage{algpseudocode}
\usepackage{caption}
\usepackage{color}
\usepackage{bbm}
\usepackage{tikz}
\usetikzlibrary{arrows.meta}
\usepackage{mathtools}
\mathtoolsset{showonlyrefs,showmanualtags}

\DeclareMathAlphabet{\pazocal}{OMS}{zplm}{m}{n}

\newcommand{\ve}{\boldsymbol}

\renewcommand{\leq}{\leqslant}

\renewcommand{\epsilon}{\varepsilon}

\newtheorem{defn}{Definition}

\newtheorem{prep}{Proposition}

\usepackage{pgfplots}
\DeclareUnicodeCharacter{2212}{−}
\usepgfplotslibrary{groupplots,dateplot}
\usetikzlibrary{patterns,shapes.arrows}

\modulolinenumbers[5]

\usepackage{multirow}
\usepackage{color}
\usepackage{paralist}
\usepackage[space]{cite}

\newtheorem{proc}{Procedure}

\pgfplotsset{compat=1.18}

\IEEEoverridecommandlockouts                        

\title{\LARGE \bf
Data-Driven Estimation of Vinnicombe metric}

\author{Margarita A. Guerrero, Henrik Sandberg,
 and Cristian R. Rojas
\thanks{}% <-this % stops a space
\thanks{This work was partly supported by the Swedish Research Council under contract number 2023-05170, by  the Wallenberg AI, Autonomous Systems and Software Program (WASP), and by the Swedish Civil Defence and Resilience Agency (Project MAD-VAMCHS). The authors are with the Division of Decision and Control Systems, KTH Royal Institute of Technology, 100 44 Stockholm, Sweden (e-mails: mags3@kth.se, hsan@kth.se, crro@kth.se). 
}
\thanks{Code available at \href{https://github.com/mags-ono/Vinnicombe_Metric}{https://github.com/mags-ono/nu-gap}.}
}

\begin{document}

\maketitle

%%%%%%%%%%%%%%%%%%%%%%%%%%%%%%%%%%%%%%%%%%%%%%%%%%%%%%%%%%%%%%%%%%%%%%%%%%%%%%%%
\begin{abstract}
Quantifying model mismatch in a control-relevant manner is fundamental in robust control. A well-known metric for this purpose is the $\nu$-gap, or Vinnicombe metric, which measures the discrepancy between a nominal model and the real system from a closed-loop viewpoint. However, its computation typically requires explicit knowledge of the true system. In this letter, we propose an identification-free, data-driven method to estimate the $\nu$-gap between discrete-time SISO systems directly from input-output experiments. The method is complemented by a data-driven winding-number test, based on Welch-type averaging, to verify a required topological condition for the computation of the metric. Numerical simulations on heavy-duty gas-turbine models and a textbook example show that the proposed estimate closely matches MATLAB$^\copyright$ \texttt{gapmetric}, while correctly detecting cases in which the admissibility conditions fail.
\end{abstract}
%%%%%%%%%%%%%%%%%%%%%%%%%%%%%%%%%%%%%%%%%%%%%%%%%%%%%%%%%%%%%%%%%%%%%%%%%%%%%%%%
\section{Introduction}
A core aim of automatic control is to design robust controllers that operate reliably despite modeling errors. In this context, robustness means that performance does not deteriorate significantly when the real system differs from the model used for the controller synthesis. Several frameworks have been developed to address this problem, including small-gain arguments~\cite{desoer_09}, integral quadratic constraints (IQCs)~\cite{Megretski_97}, and $\mathcal{H}_\infty\!$ control~\cite{zhoudoyle}. Essentially, these approaches provide stability or performance guarantees for all plants whose uncertainty satisfies a prescribed description, such as norm bounds, structural assumptions, or multiplier-based conditions.

A control-relevant way to quantify model mismatch is provided by the Vinnicombe metric~\cite{vinnicombe_91,vinnicombe_92}, also known as the $\nu$-gap metric, defined for linear time-invariant systems (LTI). This metric is particularly appealing because it measures discrepancies between plants from a closed-loop viewpoint and induces the weakest topology in which closed-loop stability is a robust property. Thus, if the $\nu$-gap between a nominal model and the true plant is sufficiently small, then a controller that robustly stabilizes the nominal model will also stabilize the true plant, under suitable admissibility conditions. Moreover, while topologically equivalent to the gap metric~\cite{sakkary_85,cantoni_17}, the $\nu$-gap admits a direct frequency-response interpretation.

%: if the $\nu$-gap between a nominal model and the true plant is sufficiently small, then a controller that robustly stabilizes the nominal model will also stabilize the true plant, under suitable admissibility conditions. 

% Yet computing the $\nu$-gap typically requires explicit knowledge of the true system, which is precisely the information that is unavailable in many practical settings.
In practice, however, the usefulness of both robust-control guarantees and $\nu$-gap-based reasoning depends on how accurately the plant uncertainty or model mismatch can be characterized. When the true plant is unknown and prior knowledge is inaccurate, constructing a meaningful uncertainty description may be difficult, conservative, or infeasible. 

% This particularly problem has encouraged researchers on finding robust metrics in data-driven settings, where only input-output measurements are available and no explicit prior or plant model is assumed~\cite{koch2022determining,koenings2017data}. One popular approach involves estimating the $\mathcal{H}_\infty$ norm of the modeling error directly from experiments on the system by employing iterative methods, such as power iterations~\cite{wahlberg2010non,rojas2012analyzing} and Thompson sampling~\cite{muller2019gain}, without the need for a parametric model. Recently a similar power-iterations approach to compute the structured singular value has been proposed in~\cite{guerrero_25}.
This practical difficulty has motivated research on estimating control-relevant robustness quantities directly from data, in settings where only input-output measurements are available and no explicit plant knowledge is assumed~\cite{koch2022determining,koenings2017data}. One prominent line of work estimates the $\mathcal{H}_\infty$ norm of the modeling error directly from experiments by means of iterative procedures, such as power iterations~\cite{wahlberg2010non,rojas2012analyzing} and Thompson sampling~\cite{muller2019gain}. %, without requiring a parametric model. 
More recently, a power-iteration-based approach for estimating the structured singular value has been proposed in~\cite{guerrero_25}.

To the best of our knowledge, there is no data-driven method for estimating the Vinnicombe metric. Motivated by this, we propose an identification-free, iterative approach that estimates the 
$\nu$-gap between two plants directly from time-domain input-output data. The proposed method builds on the core idea of power iterations to compute a numerical estimate of the distance between two dynamical systems, while an iterative winding-number calculation is used to verify the admissibility conditions under which the 
$\nu$-gap is well defined.

Our main contributions are as follows:
\begin{itemize}
    \item We propose a novel identification-free, data-driven method for estimating the Vinnicombe metric directly from time-domain input-output experiments.
    
\item We introduce an identification-free procedure for verifying the admissibility conditions of the $\nu$-gap through an iterative winding-number computation based on Welch-inspired frequency-domain estimates.
    
    \item We demonstrate the effectiveness of the proposed method through numerical simulations, showing that the resulting estimates closely agree with MATLAB\textsuperscript{\textregistered} \texttt{gapmetric}.
    
    \item We illustrate the practical applicability of the proposed framework to General Electric heavy-duty gas turbines, showing how the Vinnicombe metric can be interpreted and implemented in an industrial setting.

\end{itemize}
The remainder of this paper is organized as follows: In Section~\ref{sec: setup}, we state our problem setup. Section~\ref{sec: prelims} reviews the Vinnicombe metric and the power iteration method. In Section~\ref{sec: approach}, we outline our proposed approach, and demonstrate its efficacy in Section~\ref{sec: Simulation}. The paper is concluded in Section~\ref{sec: conclusion}.

\section{Problem Statement} \label{sec: setup}
Consider a stable unity-feedback interconnection of a nominal plant $G_0$ with a controller $C$, denoted by $[G_0,C]$. Suppose that the nominal plant differs from an unknown true system $G$. Our objective is to guarantee that the controller $C$ will perform as intended on the true process regardless of the unmodeled dynamics. Throughout the paper, $(G_0,G,C)$ are discrete-time, linear time-invariant, single-input single-output systems.%In our setup $(G_0,G,C)$ are discrete-time, linear time-invariant, single-input single-output systems.

The aforementioned problem can be handled with a well-known distance metric, namely the $\nu$-gap, also called the Vinnicombe metric~\cite{vinnicombe_91,vinnicombe_92}. The Vinnicombe metric, denoted by $\delta_\nu$, is an appropriate measure for assessing how two systems differ in their closed-loop behavior, satisfying $0\leq \delta_\nu(G_1,G_2)\leq 1$ (see Section~\ref{sec:Vinnicombe}). 
% Its formal definition is postponed to Section~\ref{sec:Vinnicombe}. 

This quantity admits a closed-loop interpretation: it measures, in a normalized worst-case sense, the discrepancy between the complementary sensitivity behaviors induced by unity feedback around $G$ and $G_0$. In particular, taking a sufficiently small value of $\delta_\nu$ implies that any controller that stabilizes $G_0$ with adequate margin also stabilizes $G$ (under additional conditions summarized in Section~\ref{sec:Vinnicombe}).

% The Vinnicombe metric is an appropriate measure for assessing how two systems differ in their closed loop behavior. 
% In particular, let us consider two systems with transfer functions $G_1(z)$ and $G_2(z)$, and define

% \begin{IEEEeqnarray}{rCl}
% \label{eq:d_v}
% d(G_1,G_2)
% := \sup_{\omega}
% \frac{\left|G_1(e^{j\omega})-G_2(e^{j\omega})\right|}
% {\sqrt{(1+\left|G_1(e^{j\omega})\right|^2)(1+\left|G_2(e^{j\omega})|^2\right)}},\;
% \end{IEEEeqnarray}

 %which is a metric satisfying $0\leq d(G_1,G_2)\leq 1$. 
 % This quantity admits a closed-loop interpretation: it measures, in a normalized
 % worst-case sense, the discrepancy between the complementary sensitivity behaviors induced by unity feedback around
 % $G_1$ and $G_2$. 
 
 % The number $d(G_1,G_2)$ can be interpreted as the difference between the complementary sensitivy functions for the closed loop systems that are obtained with unit feedback around $G_1$ and $G_2$. 
 
% In particular, taking $G_1=G_0$ and $G_2=G$, a sufficiently small value of~\eqref{eq:d_v} implies that any controller that stabilizes $G_0$ with adequate margin also stabilizes $G$ (under additional conditions summarized in Section~\ref{sec:Vinnicombe}).
% In particular, taking $G_1=G_0$ and $G_2=G$, if \eqref{eq:d_v} is sufficiently small, then $C$ is guaranteed to stabilize $G$ (along with other conditions which are covered in Section~\ref{sec:Vinnicombe}).

When the real plant is unknown, computing $\delta_\nu(G,G_0)$ becomes untenable, as it requires knowledge of the true process. A classical workaround is to postulate that $G$ lies in a neighborhood of the nominal model, i.e., $\mathcal{B}_\tau(G_0)\coloneq\{\tilde G\colon \delta_\nu(\tilde G,G_0) \le \tau\}$~\cite[Proposition~1.2]{vinnicombe_92}. However, specifying such a set
typically relies on prior modeling knowledge about unmodeled dynamics and may be conservative or inaccurate.%Such a set is usually specified using prior modeling insight on unmodeled dynamics and uncertainty structure; in practice, the resulting bound may be conservative or inaccurate.
% belongs to a $\nu$-gap ball around the nominal model, i.e.,
% $\{\,\tilde G:\delta_\nu(\tilde G,G_0)\le \tau\,\}$~\cite[Proposition~1.2]{vinnicombe_92}. However, specifying such a set
% typically relies on prior modeling knowledge about unmodeled dynamics and may be conservative or inaccurate.%This problematic could be addressed if we assume that the true system belongs to a set $\{\tilde{G}: \delta_\nu(\tilde{G},G_0)\leq \tau\}$ \cite[Preposition 1.2]{vinnicombe_92}. This last assumption require prior knowledge of the real plant which may be inaccurate. 

Motivated by these limitations, we propose in Section~\ref{sec: approach} an identification-free, data-driven procedure
to compute the Vinnicombe metric for the single-input single-output (SISO) case using a power-iteration principle~\cite{golub2013matrix}, requiring only
input-output experiments on the black-box system $G$ and an iterative redesign of the excitation signal. While we
primarily assume that $G_0$ is an available nominal model, the same methodology can be used as a data-driven similarity
test when both systems are unknown, e.g., to compare the opening/closing dynamics of two valves from closed-loop
measurements.
% Motivated by these limitations,
% we introduce in Section~\ref{sec: approach} a fully data-driven computation of the Vinnicombe metric based on the power iterations method~\cite{golub2013matrix} requiring only simulations on the black-box model $G$ and an iterative redesign of the input signal. In our setting we assume $G_0$ is a model of the real plant, but our implementation can be additionally extended for when both systems $(G,G_0)$ are unknown and a similarity test on their closed-loop response is needed, e.g. describe how two valves differs in their aperture-close response.

Before presenting the proposed approach, we introduce the technical preliminaries required throughout the manuscript.
% Before we propose our approach, we shall present several technical preliminaries that are required for this paper. 
% habalr de v-ap y su implementacion y por qué es usperior a otras metricas

% pero despues decir.. y qué pasa cuando no se cuenta con la medida real de la planta?
% It turns out that, as commonly encountered in system identification endeavor, our nominal plant differs from the unknown-true system dynamics. 

\section{Preliminaries} \label{sec: prelims}
% In this section, we review the model-based power method for computing a lower bound on $\mu_\Delta$. The results stated in this section appear in~\cite{packard1988power,doyle1982analysis}. 
\subsection{The Vinnicombe Metric}\label{sec:Vinnicombe}
\vspace{1em}
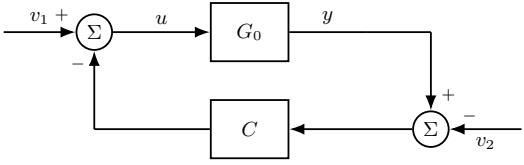
\begin{figure}[t]
\centering
\resizebox{0.8\columnwidth}{!}{%
\begin{tikzpicture}[auto, >=Latex, thick, font=\small]

    % Nodes
    \node[circle, draw, minimum size=5.5mm, inner sep=0pt] (sum1) at (0,0) {$\Sigma$};
    \node[draw, minimum width=1.2cm, minimum height=0.9cm] (P) at (2.4,0) {$G_0$};
    \node[circle, draw, minimum size=5.5mm, inner sep=0pt] (sum2) at (5.2,-1.5) {$\Sigma$};
    \node[draw, minimum width=1.2cm, minimum height=0.9cm] (C) at (2.4,-1.5) {$C$};

    % Input v1
    \draw[->] (-1.4,0) -- node[above] {$v_1$} (sum1);

    % Sum1 to plant
    %\draw (sum1) -- node[above] {$u$} (P);
    \draw[->] (sum1) -- node[above] {$u$} (P);

    % Plant output y
    \draw (P) -- node[above] {$y$} (4.2,0);
    \draw (4.2,0) -- (5.2,0);
    \draw[->] (5.2,0) -- (sum2);

    % Input v2
    \draw[->] (6.6,-1.5) -- node[below] {$v_2$} (sum2);

    % Sum2 to controller
    %\draw (sum2) -- (C);
    \draw[->] (sum2) -- (C);

    % Controller back to sum1
    \draw (C) -- (0,-1.5);
    \draw[->] (0,-1.5) -- (sum1);

    % Minus signs
    \node at (-0.25,-0.5) {\scriptsize $-$};
    \node at (-0.5,0.26) {\scriptsize $+$};

    \node at (5.8,-1.3) {\scriptsize $-$};
    \node at (5.47,-0.97) {\scriptsize $+$};

\end{tikzpicture}%
}
\vspace{1.5mm}
\caption{Standard feedback configuration.}
    \label{fig:fb} \vspace{-1.7em}
\end{figure}
% Given a nominal discrete time, linear time-invariant system  $G_0(z)$, and a feedback compensator $C(z)$, we define its transfer function from $v_1$ and $v_2$, to $u$ and $y$ (see Figure~\ref{fig:fb}) as follows
% \[\begin{bmatrix} G_0 \\ I \end{bmatrix}(I-CG_0)^{-1}\begin{bmatrix} -C & I \end{bmatrix}\]
Consider the standard negative-feedback interconnection of a nominal discrete-time LTI multi-input multi-output (MIMO) plant $G_0(z)$
and a compensator $C(z)$ as shown in Fig.~\ref{fig:fb}. Let $v_1$ and $v_2$ denote additive disturbances at
the controller input and output, respectively, and let $(y,u)$ be the plant output and input. The
closed-loop transfer matrix from $[v_1\; v_2]^T$ to $[y\; u]^T$ is
\vspace{-1em}
\begin{IEEEeqnarray}{rCl}
\label{eq:T0}
T(G_0,C)
\;:=\;
\begin{bmatrix} G_0 \\ I \end{bmatrix}
\bigl(I + C G_0\bigr)^{-1}
\begin{bmatrix} -C & I \end{bmatrix},
\end{IEEEeqnarray}
where $I=1$ is the SISO case.

We define the \emph{stability margin} of the feedback interconnection $[G_0,C]$ as
\begin{IEEEeqnarray}{rCl}
% \label{eq:Tgc}
% T(G_0,C) &:=&
% \left[\begin{array}{c} G_0 \\ I \end{array}\right](I-CG_0)^{-1}
% \left[\begin{array}{cc}-C & I\end{array}\right],
% \\[2pt]
\label{eq:sm}
b_{G_0,C} &:=&
\begin{cases}
\|T(G_0,C)\|_\infty^{-1}, & \text{if } [G_0,C]\ \text{is stable},\\
0, & \text{otherwise.}
\end{cases}
\end{IEEEeqnarray}
The number $b_{G_0,C}$ can be regarded as a measure of the performance of the feedback system comprising $G_0$ and $C$, with larger values of \eqref{eq:sm} corresponding to better performance. Note that \eqref{eq:sm} lies in the range $[0,1]$~\cite{vinnicombe_91}.

\begin{defn}
\label{def:defnu}
Let $G_0(z)$ and $G(z)$ be discrete-time single-input single-output (SISO) transfer functions. The Vinnicombe metric (or $\nu$-gap)~\cite{vinnicombe_92}
between $G_0$ and $G$ is defined as
\begin{IEEEeqnarray}{l}
\hspace*{-1em}\delta_\nu(G_0,G)= \nonumber\\[-2pt]
{
\begin{cases}
\displaystyle 
\left\|
\frac{G_0-G}{\sqrt{(1+|G_0|^2)(1+|G|^2)}}
\right\|_\infty,
& \text{if } (G_0,G)\in\mathcal{C},\\[2pt]
1, & \text{otherwise.}
\end{cases}
}
\label{eq:nugap}
\end{IEEEeqnarray}
% where $(G_0,G)\!\in \mathcal{C}$ if the index functions $f_1=1+G_0(z)G_0(z^{-1})$ and $f_2=1+G(z)G_0(z^{-1})$ have the same number of zeros in $\mathbb{E}\coloneq\{z\in \mathbb{C}: |z|>1\}$.
where $(G_0,G)\in\mathcal{C}$ if the functions
$f_1(z)\coloneq 1+G_0(z)G_0(z^{-1})$ and $f_2(z)\coloneq 1+G(z)G_0(z^{-1})$
have the same number of zeros in $\mathbb{E}\coloneq\{z\in\mathbb{C}\colon |z|>1\}$.
\end{defn}
The condition $(G_0,G)\in\mathcal{C}$ is an index (winding-number) requirement that ensures a
meaningful notion of continuity: if $G$ is perturbed continuously away from $G_0$, it rules out
situations in which an intermediate plant $\tilde G$ would necessarily yield $\delta_\nu(\tilde G,G_0)=1$
even though $\delta_\nu(G,G_0)$ is small~\cite[Section 12.1]
{aastrom_21feedback}. In the SISO case, this condition can be checked via a
Nyquist/winding-number test, where the winding number is denoted by $\mathrm{wno}(\cdot)$; in particular, the Nyquist plot of $G(e^{j\omega})G_0(e^{-j\omega})$ should
not encircle the critical point $-1$ for $\omega\in(-\pi,\pi]$. Additionally, in this same SISO setting, the pointwise quantity underlying the $\nu$-gap admits a geometric interpretation as the chordal distance between the projected Nyquist diagrams of two plants onto the Riemann sphere~\cite[Sec.~3]{vinnicombe_92}. Therefore, $\delta_\nu$ may be viewed as the maximum normalized separation over frequency, which explains its relevance as a closed-loop similarity metric.

\begin{prep}[\hspace{-0.2mm}{\cite[Prop.~1.2]{vinnicombe_92}}]
\label{prep:1} Given a nominal plant $G_0$ and a perturbed plant $G$, the interconnection $[G,C]$ is stable for all compensators $C$ satisfying $b_{G_0,C}>\beta$ if and only if $\delta_\nu(G_0,G)\leq\beta$.
% Given a nominal plant $G_0$, and a perturbed plant $G$, any controller $C$ that stabilizes $G_0$, and achieves $b_{G_0,C}>\beta$, stabilizes $[G,C]$ if and only if $\delta_\nu(G_0,G)\leq\beta$.
\end{prep}

The $\nu$-gap metric $\delta_\nu(G_0,G)$ quantifies the distance between two plants from a closed-loop
perspective. In particular, $\delta_\nu(G_0,G)\in[0,1]$, with $\delta_\nu(G_0,G)=0$ if and only if
$G_0=G$. Hence, when $G$ represents the (unknown) true plant, a small value of $\delta_\nu(G_0,G)$ indicates
that the mismatch between $G_0$ and $G$ is negligible in terms of closed-loop stability.
%Proposition~\ref{prep:1} shows that $\delta_\nu$ is precisely the metric that matches the
%stability margin: if a controller $C$ achieves a margin $b_{G_0,C}>\beta$ on the nominal plant $G_0$,
%then $C$ stabilizes any plant $G$ satisfying $\delta_\nu(G_0,G)\le \beta$, and this bound is tight.

\subsection{Power Iteration Method}
In order to compute~\eqref{eq:nugap}, one needs to compute the $\mathcal{H}_\infty$ norm of a transfer function. To this end, consider a linear discrete-time stable causal and proper system of the form
\[y_k=G(q)u_k+e_k\]
where $\{u_k\}$ is an input signal, $G$ is a transfer function, $\{y_k \}$ is the output of the system, and $\{ e_k\}$ is Gaussian, not necessarily white noise. Next, recall that the $\mathcal{H}_\infty$ norm is an induced norm, i.e., $\left\|G\right\|_\infty=\left\|G\right\|_{i2}$; see \cite[App. A.5.7]{skogestad2005multivariable} for a proof. Assume that, at iteration $t$, the signals have a finite length $N \in \mathbb{N}$, i.e., $\ve{u}^t:=[u_0^t,\dots,u_{N-1}^t]^T$, $\ve{y}^t:=[y_0^t,\dots,y_{N-1}^t]^T \in \mathbb{R}^{N\times1}$, hence
\vspace{-0.4em}
\begin{IEEEeqnarray}{rCl}\label{eq:i2N} 
\|{G}\|_{i2,N}
= \sup_{\ve{u}^t\in\mathbb{R}^{N\times 1}\setminus\{\ve{0}\}}
\frac{\|\ve{y}^t\|_2}{\|\ve{u}^t\|_2},
\end{IEEEeqnarray}
where $\|G\|_{i2,N}\rightarrow\|G\|_\infty$ for $N\rightarrow\infty$ \cite[Th. 3]{rojas2012analyzing}. Let us assume that, at iteration $t$, the system starts from  zero initial conditions. If there is no noise affecting the output of the system, the mapping from $\ve{u}^t$ to $\ve{y}^t$
%\{u_k^t\}_{k=1,\dots,N}$ to \{y_k^t\}_{k=1,\dots,N}$
can be written as
\vspace{-0.7em}
\[\ve{y}^t =\ve{G}_{N}\ve{u}^t, \vspace{-0.7em}\]
where %$\ve{u}=[u_1^t \cdots u_N^t]^T$, $\ve{y}=[y_1 \cdots y_N^t]^T$ and
\vspace{-0.7em}
\[
\ve{G}_N=
\begin{bmatrix}
g_0     & 0      & 0      & \cdots & 0 \\
g_1     & g_0    & 0      & \cdots & 0 \\
\vdots  & \ddots & \ddots & \ddots & \vdots \\
g_{N-1} & g_{N-2}& \cdots & g_1    & g_0
\end{bmatrix}
\]
is a triangular Toeplitz matrix. Thus, \eqref{eq:i2N} is equivalent to
\begin{IEEEeqnarray}{rCl} \|G\|_{i2,N}^2
 \approx \sup_{\ve{u}^t \neq0}
\frac{(\ve{u}^t)^{T}\ve{G}_N^{T}\ve{G}_N\,\ve{u}^t}{(\ve{u}^t)^{T}\ve{u}^t}
= \lambda^{\max}\!\left(\ve{G}_N^{T}\ve{G}_N\right).\end{IEEEeqnarray}
The right hand side, known as the \textit{Rayleigh quotient} of $\ve{G}^T\ve{G}$, corresponds to the largest eigenvalue of $\ve{G}^T\ve{G}$, and it can be computed using the power iteration method.
% Thus, for $N$ sufficiently large,
% \begin{IEEEeqnarray}{rCl}
% \label{eq:inftynorm}\|\ve{G}\|_{i2}^2  \approx \sup_{\ve{u} \neq0}
% \frac{\ve{u}^{T}\|G\|^2\,\ve{u}}{\ve{u}^{T}\ve{u}}
% = \lambda_{\max}\!\left(\|G\|^2\right).\end{IEEEeqnarray}
% % \begin{IEEEeqnarray}{rCl}
% % \label{eq:inftynorm}\|\ve{G}\|_{i2} \xrightarrow{N\to\infty}
% % \|\ve{G}\|_{\infty} :=\sup_{\omega\in (-\pi,\pi]}
% % |G(e^{j\omega})|\end{IEEEeqnarray}
% Our model-free framework requires \eqref{eq:inftynorm} to be computed without a parametric model. One could appeal to the so-called \textit{power iterations method}~\cite{rojas2012analyzing} from linear algebra.% to estimate the $\mathcal{H}_\infty$-norm of $G$.
\begin{proc}[Power Iteration Method]\label{pro:pm}
Consider a square matrix $\ve{A}\in \mathbb{R}^{n\times n}$. The power iteration method performs the following steps to estimate the eigenvector of $\ve{A}$ corresponding to its largest eigenvalue:
\begin{itemize}
    \item[1.] Choose $\ve{x}_0 \in \mathbb{R}^n$ randomly.
    \item[2.] Set $i \gets 1$.
    \item[3.] Compute $\ve{x}_i \gets \frac{1}{\|\ve{A}\ve{x}_{i-1}\|_2}\ve{A}\ve{x}_{i-1}$.
    \item[4.] Set $i\gets i+1$ and go to Step 3.
\end{itemize}

For normalized iterates, a natural estimate of the magnitude of the dominant eigenvalue is $\|\ve{A}\ve{x}_{i-1}\|_2$\footnote{If $\ve{A}$ is symmetric, the dominant eigenvalue may also be estimated through the Rayleigh quotient $\ve{x}_{i-1}^T\ve{A}\ve{x}_{i-1}$. The two estimates may differ numerically under noise or finite-iteration effects.}.
% The largest eigenvalue of $\ve{A}$, at iteration $i$, can be estimated as $\hat{\lambda}^{\max}=\ve{x}_{i-1}^T\ve{A}\ve{x}_{i-1}$.
\end{proc}
In the SISO frequency-domain setting, the $\mathcal{H}_\infty$ norm coincides with the maximum magnitude of the frequency response over all frequencies, i.e., with the peak of the Bode magnitude plot. %Since $G(e^{j\omega})$ is scalar-valued, its largest singular value is simply its magnitude.
Thus, the same power-iteration principle can be interpreted in the frequency domain as an iterative search
for the peak gain, where the maximizing frequency is implicitly identified by the iterations.

% In order to close the gap between the power iteration in Procedure~\ref{pro:pm} with \eqref{eq:inftynorm} we recall that for the SISO case
% \[\|G\|_\infty = \]
Our proposed procedure retains this template to approximate the $\nu$-gap in~\eqref{eq:nugap} based solely on experiments on the unknown system $G$. We assume $e_k = 0$ throughout the analysis, while noise is included in the simulations of Section~\ref{sec: Simulation} to evaluate robustness.

\section{Proposed Approach} \label{sec: approach}
In this section, we propose a data-driven procedure to estimate the Vinnicombe $\nu$-gap metric from time-domain experiments. Assuming that a nominal model $G_0$ is available and that only input-output data can be collected, we adapt the model-based construction in~\cite{vinnicombe_91} to operate directly on frequency-domain quantities estimated from time-domain data.
\vspace{-2.2em}
\subsection{$\nu$-gap estimation}

We focus on the practically relevant setting where $G$ denotes the (unknown) plant and $G_0$ is a nominal model. For a given input $u$, let $U(e^{j\omega})$ be its discrete-time Fourier transform~(DFT), and let $Y(e^{j\omega})$ and $Y_0(e^{j\omega})$ denote the corresponding output spectra obtained from experiments on $G$ and $G_0$, respectively, so that $Y(e^{j\omega})=G(e^{j\omega})U(e^{j\omega})$ and $Y_0(e^{j\omega})=G_0(e^{j\omega})U(e^{j\omega})$. When $(G,G_0)\in\mathcal{C}$, the $\nu$-gap admits the equivalent frequency-domain expression
\begin{IEEEeqnarray*}{l}
\delta_{\nu}(G,G_0)
=\\
\sup_{\omega \in (-\pi,\pi]}
\frac{|U(e^{j\omega})|^2\,|Y(e^{j\omega})-Y_0(e^{j\omega})|}
{\sqrt{|U(e^{j\omega})|^2+|Y(e^{j\omega})|^2}\;
 \sqrt{|U(e^{j\omega})|^2+|Y_0(e^{j\omega})|^2}}.
\end{IEEEeqnarray*}
As in $\mathcal{H}_\infty$ norm estimation, the $\nu$-gap is the supremum of the modulus of a frequency function. In the noiseless case, this motivates a power-iteration-type input redesign that progressively concentrates excitation at frequencies that dominate the above ratio. Specifically, given spectra $(U^n,Y^n,Y_0^n)$ at iteration $n$, we update the next input in the frequency domain according to
\begin{IEEEeqnarray}{l}
\tilde{U}^{n+1}\!\left(e^{j\omega}\right)= \nonumber\\
\frac{\left|U^{n}\!\left(e^{j\omega}\right)\right|^{2}\left[\,Y^{n}\!\left(e^{j\omega}\right)-Y_{0}^{n}\!\left(e^{j\omega}\right)\right]}
{\sqrt{\left|U^{n}\!\left(e^{j\omega}\right)\right|^{2}
+\left|Y^{n}\!\left(e^{j\omega}\right)\right|^{2}}\;
 \sqrt{\left|U^{n}\!\left(e^{j\omega}\right)\right|^{2}
+\left|Y_{0}^{n}\!\left(e^{j\omega}\right)\right|^{2}}}\!.
\label{eq:tildeu}
\end{IEEEeqnarray}
The role of \eqref{eq:tildeu} is then analogous to that of the matrix-vector product in the classical power iteration: each update amplifies the frequency components that contribute most to the normalized discrepancy between $G$ and $G_0$. The updated spectrum is mapped back to the time domain via an inverse DFT, normalized, and used in the next experiment. Thus, the input progressively concentrates near the frequency at which the $\nu$-gap is attained.
\vspace{-1em}
\subsection{$\mathrm{wno}$ estimation}
The update rule in \eqref{eq:tildeu} suggests a straightforward implementation; however, the index condition $(G_0,G)\in\mathcal{C}$ cannot be overlooked. This condition rules out pathological situations where, despite a small closed-loop mismatch between $G_0$ and $G$, any continuous path from $G_0$ to $G$ would necessarily cross an intermediate transfer function $\widetilde G$ for which the $\nu$-gap saturates to one. Accordingly, we verify that $(G_0,G)\in\mathcal{C}$ in a data-driven fashion by computing winding numbers of the index functions $f_1$ and $f_2$, i.e., by checking that the Nyquist curve of $G(e^{j\omega})G_0(e^{-j\omega})$ does not encircle the critical point $-1$ for $\omega\in(-\pi,\pi]$.
%Even thought the implementation from now seems straightforward, the assumption that $(G_0,G)\in \mathcal{C}$, i.e., that there is no intermediate transfer function $\tilde{G}$ such that $\delta_\nu (G_0,\tilde{G})$ is $1$ even if $\sigma_\nu (G_0,G)$ is small, no se puede pasar por alto, dado que esta es la limitación más importante cuando se calcula el vinnicombe metric con la finalidad de comparar el comportamiento de los sistemas under feedback. Accordingly, we calculate in a data-driven fashion the winding-number on $f_2$ (and $f_1$) to verify $G(e^{j\omega})G_0(e^{-j\omega})$ do not encircle $-1$ for $w\in(-\pi,\pi]$.

Using the same information from the Fourier transform of the input and outputs of~\eqref{eq:tildeu}, and inspired by Welch's method~\cite{welch_67} we compute $f_2$ as follows:
\begin{IEEEeqnarray}{l}
f_2\!\left(e^{j\omega_k}\right) \approx \nonumber\\ 1 +
\frac{
\left(\sum_{n=1}^{N_{\mathrm{acc}}} Y_{n}[k]\,U_n[k]^{*}\right)
\left(\sum_{n=1}^{N_{\mathrm{acc}}} Y_{0,n}[k]\,U_n[k]^{*}\right)^{*}
}{
\left(\max\!\left(\sum_{n=1}^{N_{\mathrm{acc}}}\lvert U_n[k]\rvert^{2},\,\varepsilon_{0}\right)\right)^{2}
},
\label{eq:f2}
\end{IEEEeqnarray}
where $n$ is the current iteration,  $\omega_k:=2\pi k/N$, $k\in\{0,1,\dots,N-1\}$ is the DFT grid on the unit circle, $N_\text{acc}$ is the number of iterations used to approximate $f_2$, and $\epsilon_0$ is a small number to avoid division by 0. An analogous expression can be written for $f_1$ by replacing $(Y_n,Y_{0,n})$  with $(Y_{0,n},Y_{0,n})$ in~\eqref{eq:f2}. Inspired by Welch’s averaged periodogram estimator~\cite[Ch.~6, Eq.~(6.70)]{Ljung:99}, we estimate the index functions by averaging across power-iteration experiments. Unlike classical Welch's method, which partitions a single record into multiple sub-batches and recomputes the DFT for each batch, our approach reuses the DFTs that are already computed at each power-iteration step. We therefore average over the first $N_{\mathrm{acc}}\ll M$ iterations: as the power method progresses, the redesigned input becomes increasingly narrowband around the peak frequency, which reduces spectral coverage and makes later iterations less informative for index verification~\cite{rojas2012analyzing}.
% The main difference between the batches approach of the Welch method~\cite[Ch. 6, Eq. (6.70)]{Ljung:99} is that in \eqref{eq:f2} we average over the power-method iterations, reusing the data for the iterative method and avoiding computing the Discrete Fourier transform again for each batch. If the total number of iterations is $M$, then $N_\textbf{acc}\ll M$ given by the fact that over iterations the expectrum of the output redesigned signal focus on frequencies that contribute most to the maximum singular value, i.e. the peak frequency, losing information of rest of the spectrum. 
% Finally, we compute the winding numbers $\mathrm{wno}(f_2)$ and $\mathrm{wno}(f_1)$ by summing discrete phase increments along the closed contour and dividing by $2\pi$. Specifically, letting $f_i[N]=f_i[0]$ to close the contour, we define
%
% %otra opcion
% Finally, letting $f_i[k]:=f_i(e^{j\omega_k})$, with $\omega_k=2\pi k/N$ for $k=0,\dots,N-1$,
% Finally, letting $f_i[k]:=f_i(e^{j\omega_k})$ denote the samples of $f_i$ on the DFT grid, we compute the winding numbers $\mathrm{wno}(f_2)$ and $\mathrm{wno}(f_1)$ by summing discrete phase increments along the closed contour obtained from the sampled points $\{f_i[k]\}_{k=0}^{N-1}$, where we set $f_i[N]:=f_i[0]$ to close the contour. Specifically, we define
Finally, for $i\in\{1,2\}$, let $f_i[k]:=f_i(e^{j\omega_k})$ denote the samples of $f_i$ on the DFT grid. We compute the winding numbers $\mathrm{wno}(f_2)$ and $\mathrm{wno}(f_1)$ by summing discrete phase increments of $f_i$ around the unit circle, where we set $f_i[N]:=f_i[0]$ to close the sampled path. Specifically, we define
\begin{IEEEeqnarray}{l}
\Theta_i=\sum_{k=0}^{N-1}\arg\!\left(\frac{f_i[k+1]}{f_i[k]}\right),\;
\mathrm{wno}(f_i)=\mathrm{round}\!\left(\frac{\Theta_i}{2\pi}\right)\!.
\label{eq:f2angle}
\end{IEEEeqnarray} % and the same procedure applies to $f_1$. 
The data-driven winding-number computation is summarized in Algorithm~\ref{alg:Ccheck}.

These considerations lead to the pseudo-code for the data-driven computation of the Vinnicombe metric in Algorithm~\ref{alg:alg1}. The procedure terminates early if $\mathrm{wno}(f_1)\neq \mathrm{wno}(f_2)$ within the first $N_{\mathrm{acc}}$ iterations (indicating $(G_0,G)\notin\mathcal{C}$); otherwise, it runs for the full $M$ iterations.
% Finally we obtain the winding number 
% $wno(f_2)$ by summing the phase increments and dividing by $2\pi$, i.e.,
% \begin{IEEEeqnarray}{l}
% \Theta_i=\sum_{k=0}^{N-1}\arg\!\left(\frac{f_i[k+1]}{f_i[k]}\right),\;
% \mathrm{wno}(f_i)=\mathrm{round}\!\left(\frac{\Theta_i}{2\pi}\right)\!,
% \label{eq:f2angle}
% \end{IEEEeqnarray}
% where $i\in\{1,2\}$ as same applies for the calculation of $f_1$. The algorithm of the data-driven winding number calculation is outlined in Algorithm~2.

% The previous discussions finally lead to the pseudo-code for the data-driven computation of the Vinnicombe metric shown in Algorithm~\ref{alg:alg1}. The algorithm terminates when $wno(f_1)\neq wno(f_2)$ in the first $N_\text{acc}$ iterations or after the total number iterations $M$. 

\begin{algorithm}[h!]
\footnotesize
\caption{Data-driven estimation of the $\nu$-gap}
\label{alg:alg1}
\begin{algorithmic}[1]
\Require  A random input vector $\ve{u}^0 \in \mathbb{R}^N$
%\State Choose a random input vector $\ve{u}^0 \in \mathbb{R}^N$
\For{$n=0,1,2,\dots, M-1$}
    \State Apply $\ve{u}^n$ to $G$ and $G_0$, and collect outputs $\ve{y}^n$ and $\ve{y}_0^n$, respectively.
    \State
    Compute the DFT of $\ve{u}^n$, $\ve{y}^n$ and $\ve{y}_0^n$, and calculate~\eqref{eq:tildeu}, for each frequency $\omega$.
    \State \textbf{if} $n < N_{\mathrm{acc}}$, update accumulators via \vspace{-0.5em}
\[S_{yu}\;+\!=Y_n[k]U_n[k]^\ast, \; \;  S_{{y_0}u}\;+\!=Y_{0,n}[k]U_n[k]^\ast,  \; \; S_u\;+\!=|U_n[k]|^2\vspace{-0.5em}\]

    \textbf{end if}
    % \If{$n < N_{\mathrm{acc}}$}
    %     \State Update accumulators (batch averaging) via (X)
    % \EndIf
    \If{$n=N_{\mathrm{acc}}$}
        \State \textbf{Run} Algorithm~\ref{alg:Ccheck}
        \State \textbf{if} \texttt{inC} is false, \textbf{Terminate} and set $\delta_\nu \gets 1$ \textbf{end if}
        % \If{\texttt{inC} is false}
        %     \State \textbf{Terminate} and set $\delta_\nu \gets 1$.
        %     \State \textbf{break}
        % \EndIf
    \EndIf
    
    \State Compute $\tilde{\ve{u}}^{n+1}$ as the IDFT of $\tilde{U}^{n+1}(e^{j\omega})$
    \State Normalize $\tilde{\ve{u}}^{n+1}$ by its $2$-norm to obtain $\ve{u}^{n+1}$
    
\EndFor
\State Estimate $\delta_\nu$ as the $2$-norm of $\tilde{u}^M$
\end{algorithmic}
\end{algorithm}

\begin{algorithm}[h!]
\footnotesize
\caption{Data-driven index condition $(G_1,G_2)\in\mathcal{C}$}
\label{alg:Ccheck}
\begin{algorithmic}[1]
\Require Accumulators $(S_{yu},S_{{y_0}u},S_u)$, $N_{\mathrm{acc}}$ batches; tolerance $\varepsilon_0>0$; threshold $\texttt{tol}_f>0$.
\State \textbf{for} $n=0,1,\dots,N_\text{acc}$,  Construct $f_1$ and $f_2$ via~\eqref{eq:f2}, \textbf{end for}
%\For{$n=0,1,2,\dots, N_\text{acc}$}
%\State Construct $f_1(e^{jw_k})$ and  $f_1(e^{jw_k})$ via~\eqref{eq:f2}
%\EndFor
\State Compute minimum distances to the origin: $
m_1=\min_k |f_1[k]|$ and $ m_2=\min_k |f_2[k]|$.
\State Compute winding numbers via~\eqref{eq:f2angle}, with $f_i[N]=f_i[0],\;i\in\{1,2\}$% to close the contour, 
\State Set \texttt{inC} $\gets$ $\big(m_1>\texttt{tol}_f\big)\wedge\big(m_2>\texttt{tol}_f\big)\wedge\big(\mathrm{wno}(f_1)=\mathrm{wno}(f_2)\big)$
% \State (\textbf{Optional frequency smoothing}) Define a Hann kernel
% \[
% h[r]=\tfrac12\Big(1-\cos\frac{2\pi r}{L_s-1}\Big),\quad r=0,\dots,L_s-1,\qquad
% h \leftarrow h / \sum_{r} h[r].
% \]
% \State Compute smoothed spectra (convolution across frequency bins):
% \[
% \widetilde S_{uu}=h*S_{uu},\quad
% \widetilde S_{y_1u}=h*S_{y_1u},\quad
% \widetilde S_{y_2u}=h*S_{y_2u}.
% \]

% \State Form data-driven frequency responses (H1-type):
% \[
% \widehat H_1[k]=\frac{\widetilde S_{y_1u}[k]}{\max(\widetilde S_{uu}[k],\varepsilon_0)},\qquad
% \widehat H_2[k]=\frac{\widetilde S_{y_2u}[k]}{\max(\widetilde S_{uu}[k],\varepsilon_0)}.
% \]

% \State Construct
% \[
% f_1[k]=1+|\widehat H_1[k]|^2,\qquad
% f_2[k]=1+\widehat H_2[k]\widehat H_1[k]^*.
% \]

% \State Reorder bins to traverse $\omega\in[-\pi,\pi]$:
% \[
% f_1 \leftarrow \mathrm{fftshift}(f_1),\qquad f_2 \leftarrow \mathrm{fftshift}(f_2).
% \]

% \State Compute minimum distances to the origin:
% \[
% m_1=\min_k |f_1[k]|,\qquad m_2=\min_k |f_2[k]|.
% \]

% \State Compute winding numbers using phase increments:
% \[
% \Theta_i=\sum_{k=0}^{N-1}\arg\!\left(\frac{f_i[k+1]}{f_i[k]}\right),\qquad
% \mathrm{wind}(f_i)=\mathrm{round}\!\left(\frac{\Theta_i}{2\pi}\right),\quad i\in\{1,2\},
% \]
% with $f_i[N]=f_i[0]$ to close the contour.

% \State Set \texttt{inC} $\gets$ $\big(m_1>\texttt{tol}_f\big)\wedge\big(m_2>\texttt{tol}_f\big)\wedge\big(\mathrm{wind}(f_1)=\mathrm{wind}(f_2)\big)$.

\end{algorithmic}
\end{algorithm}

\section{Experiments} \label{sec: Simulation}
This section presents simulation studies to assess Algorithms~\ref{alg:alg1}–\ref{alg:Ccheck}. The study includes: (i) two heavy-duty gas-turbine case studies, where $\nu$-gap estimates are benchmarked against MATLAB’s \texttt{gapmetric}; and (ii) a synthetic/textbook example for which the winding number condition is analytically verifiable.

% This section presents a comprehensive set of simulations to assess the performance of Algorithm~\ref{alg:alg1} on an industrial application, specifically, a simplified Heavy-Duty Gas Turbine model.
%In this section, we present the results of extensive simulations of Algorithm~\ref{alg:power_method}, where two types of random tests are compared against the lower bound provided by \texttt{mussv} from the Robust Control Toolbox \cite[Ch.~10]{balas2012robust}. For clarity, we will refer to this lower bound as $\mu_M$ throughout the discussion.

\subsection{Case Study: Heavy-Duty Gas Turbine}
% Heavy Duty gas turbines (HDGT) are large industrial-type gas turbines used essentially for power generation by providing torque with adjustable speeds that is used to rotate a generator. HDGTs consist of an \textit{axial-flow compressor} that increases the pressure and temperature of the intake air, \textit{combustor chambers} in which fuel is injected and combusted to add thermal energy, and a \textit{turbine} section that expands the hot gases to deliver mechanical work at the shaft, typically coupled to an electrical generator.
Heavy-duty gas turbines (HDGTs) are large industrial gas turbines used mainly for power generation. They consist of an axial-flow compressor, combustor chambers where fuel is burned, and a turbine section that expands the hot gases to produce shaft work, typically driving an electrical generator.
%HDGTs consist of three primary, interconnected section: axial flow compressor, combustor chambers and turbine.
A key control variable is the Gas Control Valve (GCV), which regulates the fuel flow into the combustor and therefore directly affects the generated power. In this work, we adopt the well-known Rowen representation for heavy-duty gas turbines~\cite{rowen_83} and consider a reduced input-output model from the (measured) GCV position $u_{\mathrm{GCV}}$ to generated power $P$, obtained by cascading: (i) a fuel dynamics $G_{f,\mathrm{cl}}$, (ii) a combustor delay, (iii) a compressor-discharge dynamics $G_{cd}$, and (iv) the fuel-to-torque (power) map $F_2$.

Specifically, the fuel-path dynamics $W_f(s)$ are modeled as
\vspace{-0.7em}
\begin{equation}
W_f(s) \;=\; G_{cd}(s)\,e^{-s T_{CR}}\,G_{f,\mathrm{cl}}(s)\,u_{\mathrm{GCV}}(s),
\label{eq:wf_blocks}\vspace{-0.5em}
\end{equation}
with 
\vspace{-0.8em}
\begin{equation}
G_{f,\mathrm{cl}}(s)=\frac{1}{T_f s+1+K_F},
\qquad
G_{cd}(s)=\frac{1}{T_{cd}s+1},
\label{eq:fuel_blocks}\vspace{-0.4em}
\end{equation}

where $T_f$ is the fuel-system time constant, $K_F$ the fuel-system feedback gain, $T_{CR}$ the combustion reaction time delay, and $T_{cd}$ the compressor-discharge time constant.

The (affine) torque/power map is given by
\vspace{-0.5em}
\begin{IEEEeqnarray}{rCl}
\label{eq:f2_affine}
P \;=\; F_2(W_f,N) \;=\; A + B\,W_f + C(1-S),\vspace{-0.5em}
\end{IEEEeqnarray}

where $S$ is the per-unit rotor speed and $\{A,B,C\}$ are Rowen torque-block parameters.

For grid-connected operation, speed deviations are small and we take $S \approx 1$ over the considered operating windows. Linearizing \eqref{eq:f2_affine} around an operating point and using incremental signals (denoted by $\Delta(\cdot)$) yields the LTI incremental plant
\vspace{-1.4em}
\begin{IEEEeqnarray}{rCl}
\label{eq:Gt}
G(s)
:=
\frac{\Delta P(s)}{\Delta u_{\mathrm{GCV}}(s)}
=
\frac{B\,e^{-s T_{CR}}}{(T_f s+1+K_F)(T_{cd}s+1)}.\vspace{-0.5em}
\end{IEEEeqnarray}

% A key control variable is the Gas Control Valve (GCV), which regulates the fuel flow into the combustor and therefore directly affects the generated power. In this work, we adopt the well-known Rowen representation for heavy-duty gas turbines and derive a reduced linear time-invariant (LTI) model from the (measured) GCV position to generated power by cascading three blocks: (i) a fuel-system lag $G_{f,\mathrm{cl}}$ capturing the dominant fuel-path dynamics (optionally including fuel feedback), (ii) a compressor-discharge lag $G_{cd}$ accounting for plenum/duct dynamics, and (iii) the static fuel-to-power map $f_2$. The system can be written as follows,
% \[
% G_{\mathrm{GCV}\to P}(s)
% = f_2(W_f)= f_2(G_{cd}(s)\cdot G_{f,\mathrm{cl}}(s)),
% \]
% with 
% \begin{align}
% G_{f,\mathrm{cl}}(s)&=\frac{1}{T_f s+1+K_F},\qquad 
% G_{cd}(s)=\frac{1}{T_{cd}s+1},\\
% f_2 &= A+B\cdot W_f+C(1-N),
% \end{align}
% where $T_f$ is the fuel system time constant, $K_f$ the fuel system feedback, $T_{cd}$ the compressor discharge constant, $N$ the per unit turbine rotor speed and $\{A,B,C\}$ the gas turbine torque block parameters.
% For a turbine connected to the electrical grid $N=1$, then the resulting LTI system for an incremental plant is reduced to
% \begin{equation}
% \label{eq:Gt}
% G(s)
% :=
% \frac{\Delta P(s)}{\Delta u_{\mathrm{GCV}}(s)}
% =
% \frac{B}{(T_f s+1+K_F)(T_{cd}s+1)}.
% \end{equation}

\vspace{-1em}
\subsection*{  A.1. Nominal Model}

In~\cite{tavakoli_09}, the parameters of Rowen's model for HDGT are estimated for a 172-MW simple cycle, single shaft power unit by using
operational data at a nominal point. Following the setup in Section~\ref{sec: setup}, we treat the parameter values in~\cite[Table VII]{tavakoli_09} as the closest available description of the true plant $G$.

%a proxy for the closest available description of the true plant $G$.

% Following the idea presented in Section~\ref{sec: setup}, we assume these parameters, as shown in Table~I, to be the closest for a real system $G(z)$.

Next, we assume that a nominal transfer function $G_0$ from the (measured) GCV position $u_{\mathrm{GCV}}$ to generated power $P$ has been identified by means of a multi-objective differential evolution optimization algorithm as described in~\cite{khormali_15}. The parameters in~\cite[Table II]{khormali_15}, are estimated based on real data acquired from a gas turbine of 162 MW nominal power. 
% Let us consider now the case that we have estimated a nominal transfer function $G_0$ between the GCV position $u_{\mathrm{GCV}}$ to the generated power $P$ by means of a multi-objective differential evolution optimization algorithm as described in~\cite{khormali_15}. The parameters, as shown in Table~2, are estimated based on real data acquired from an industrial gas turbine of 162 MW nominal power. 

Assume that we have derived a controller $C$ that stabilizes $G_0$ on a high-fidelity industrial simulator. The objective is to assess whether deploying $C$ on the true operating power station is justified from a closed-loop robustness viewpoint, by estimating $\delta_\nu(G_0,G)$ from data.
To emulate realistic operating constraints, we consider that the plant is under \emph{dispatch} instructions\footnote{In power systems, ``dispatch'' refers to grid-operator instructions on power setpoints and allowable load variations.}, so that only small variations in generated power are permitted during testing.  %Consequently, only small-amplitude power perturbations are allowed during testing. 
Under these constraints, we apply a bounded excitation signal $\ve{u}^0 \in \mathbb{R}^N$ (small-signal variations around the operating point) to both $G$ and $G_0$, collect the corresponding input-output data, and run Algorithm~\ref{alg:alg1} to obtain a data-driven estimate of the $\nu$-gap and the associated index-condition check.

\begin{figure}[h!]
\centering

    %\vspace{0.4 em}
\includegraphics[width=0.75\columnwidth]{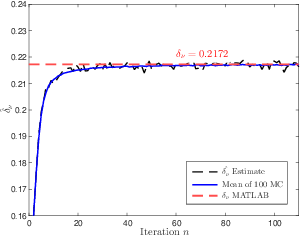}
\caption{Estimate of the $\nu$-gap $\hat{\delta}_\nu$ (black, dashed), average of 100 Monte Carlo (MC) simulations (blue, solid), and actual $\delta_\nu$ (red, dashed).}
    \label{fig:nominalA1}
\end{figure}

As shown in Fig.~\ref{fig:nominalA1}, the Monte Carlo (MC) average of the proposed estimate $\hat{\delta}_\nu$, computed over 100 runs, becomes nearly settled after roughly $60$ iterations and approaches MATLAB's reference value $\delta_\nu=0.2172$, indicating convergence of the proposed power-iteration scheme. 
The estimate is computed using $N_{\mathrm{acc}}=10$, $N=9000$, $T_s=0.05$, and measurement-noise variance $\sigma_y^2=0.01$ on both plants.

% As shown in Fig.~\ref{fig:nominalA1}, after roughly $70$ iterations the proposed estimate $\hat{\delta}_\nu$ becomes nearly settled and approaches MATLAB's reference value $\delta_\nu=0.2172$, indicating convergence of the proposed power-iteration scheme. The estimate is computed using $N_{\mathrm{acc}}=10$, $N=9000$, $T_s=0.05$, and measurement-noise variance $\sigma_y^2=0.01$ on both plants.

To further interpret this result, we inspect the dominant frequency in the last iteration and estimate the nominal local gain as $
p_0(\omega_{\mathrm{peak}})\approx \frac{Y_0(\omega_{\mathrm{peak}})}{U(\omega_{\mathrm{peak}})},$
where $Y_0$ and $U$ are the DFTs of the nominal output and input, respectively. In this case, the dominant peak occurs at $\omega_{\mathrm{peak}}\approx 0.136$ rad/sample, with $|p_0(\omega_{\mathrm{peak}})|\approx 0.54$. According to Vinnicombe's Riemann-sphere interpretation~\cite[Section 3]{vinnicombe_92}, a $\nu$-gap around $0.2$ corresponds to a relatively small discrepancy near unity gain, namely, at frequencies where $|p_0(\omega)|\approx 1$. In the present case, however, the maximizing frequency lies below unity gain, since $|p_0(\omega_{\mathrm{peak}})|\approx 0.54$. Still, the value $\hat{\delta}_\nu=0.2172$ suggests that the nominal and real plants are reasonably close from a closed-loop viewpoint. 
In particular, if the deployed controller $C$ satisfies $b_{G_0,C} > 0.2172$, then Proposition~\ref{prep:1} guarantees that it also stabilizes the true plant.

%In particular, by Proposition~\ref{prep:1}, any controller $C$ satisfying $b_{G_0,C}>0.2172$ for the nominal plant is guaranteed to stabilize the real plant as well.

% According to Vinnicombe's Riemann-sphere interpretation~\cite[Section 3]{vinnicombe_92}, a $\nu$-gap around $0.2$ corresponds to a relatively small discrepancy near unity gain ($|p_0(w_{peak})|=1$); here, although the maximizing frequency lies below unity gain, the value $\hat{\delta}_\nu=0.2172$ still suggests that the nominal and real plants are reasonably close from a closed-loop viewpoint, with a moderate mismatch concentrated in a low-frequency sub-crossover region.

% Now, as a Control System Engineer, you would like to try this controller in the real plant but, without any such performance assurances, using a designed controller directly on the physical system could potentially yield catastrophic consequences. Fortunately, the dispatch has allowed you to make changes on the power generated to realize simulations with your real system as long as these MW changes are not big in magnitude. With these guidelines you proceed to input a signal $\ve{u}^0$ to the real system and to the Nominal Model and proceed apply Algorithm~\ref{alg:alg1} accordingly.

\subsection*{  A.2. Unknown Models}
Rowen's model provides simplified dynamic representations for a broad range of General Electric (GE) heavy-duty, single-shaft gas turbines and accommodates both liquid- and gas-fuel systems. In particular, GE's Frame~6F and Frame~9F are representative mid- and high-power units (approximately 88~MW and 294~MW, respectively) for which Rowen reports ``simplified'' parameter sets and associated fuel/dynamic characteristics; see~\cite[Notes 4-7]{rowen_83}.
%
% With the rapid growth of weather-dependent solar and wind generation, gas-fired power plants are increasingly required to provide renewable-balancing services, making operational flexibility a critical requirement. GE's Dry Low NO\textsubscript{x} (DLN) combustor systems~\cite{NOX_96} support such flexible operation while meeting tightening emissions regulations, and have been widely deployed in the field. Recent upgrades of F-class DLN technology, including the DLN2.6+ platform, have further improved operability, emissions performance, and load turndown capability~\cite{venka_11}.
%
% Upgrading and maintaining DLN performance typically requires periodic tuning by a technical advisor (TA), with an on-site tuning campaign that may take on the order of 1--3 days depending on the scope (e.g., seasonal retuning due to ambient-condition changes, combustor hardware modifications, or commissioning of a new installation). For Frame~6F/9F units, the tuning of the gas control valves (GCVs) is central to DLN performance, since the GCVs regulate the fuel (and, in coordination with other actuators, the overall combustor operating condition) that ultimately determines emissions and load response.
%
With the rapid growth of weather-dependent solar and wind generation, gas-fired power plants are increasingly required to provide renewable-balancing services, making operational flexibility a critical requirement. GE's Dry Low NO\textsubscript{x} (DLN) combustor systems~\cite{NOX_96}, including recent F-class upgrades such as the DLN2.6+ platform~\cite{venka_11}, support such flexible operation while meeting tightening emissions regulations. Maintaining this performance, however, typically requires periodic tuning by a technical advisor (TA), often through on-site campaigns lasting 1--3 days depending on the scope. For Frame~6F/9F units, tuning the gas control valves (GCVs) is central to DLN performance, since these valves regulate the fuel flow and, together with other actuators, determine the combustor operating condition, emissions, and load response.

We consider a natural-gas power station equipped with both a Frame~9F and a Frame~6F unit (not necessarily within the same train). Suppose the 9F unit has been successfully tuned to comply with NO\textsubscript{x} constraints, and the operator aims to achieve comparable emissions compliance on the 6F unit. In practice, however, limited TA availability and dispatch constraints may prevent executing a full emissions test campaign for the DLN2.6+ procedure. As a result, the plant seeks an alternative, low-impact procedure to assess whether the existing GCV tuning/controller settings developed for the 9F are compatible with the 6F, while remaining within allowable small power variations approved by dispatch.

Let $G_1(z)$ and $G_2(z)$ denote the (unknown) discrete-time transfer functions from the GCV input $u_{\mathrm{GCV}}$ to generated power $P$ for the 9F and 6F units, respectively. By applying Algorithm~\ref{alg:alg1} to short, bounded tests on each unit, the maintenance team can obtain a data-driven estimate of the $\nu$-gap between $G_1$ and $G_2$ and thus quantify the degree of closed-loop mismatch. This enables an assessment of whether a GCV controller that stabilizes $G_1$ is likely to deliver comparable performance on $G_2$, i.e., whether the plant-to-plant mismatch is sufficiently small for safe deployment. Using the parameter sets in~\cite[Notes 4-7]{rowen_83}, we proceed to estimate $\delta_\nu(G_1,G_2)$. %in a fully data-driven manner.
\begin{figure}[h!]
\centering
    \vspace{-0.7 em}
    %\vspace{0.4 em}
\includegraphics[width=0.75\columnwidth]{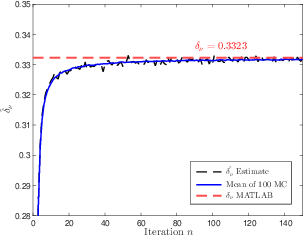}
\caption{Estimate of the $\nu$-gap $\hat{\delta}_\nu$ (black, dashed), average of 100 MC simulations (blue, solid), and actual $\delta_\nu$ (red, dashed).}
    \label{fig:unknown} 
\end{figure}

% As in the previous case, the same parameters $N_{\mathrm{acc}}=10$, $N=9000$, $T_s=0.05$, and $\sigma_y^2=0.01$ are used.
% % In Fig.~\ref{fig:unknown}, the estimate stabilizes after about $100$ iterations and converges toward MATLAB's value $\delta_\nu=0.3323$. 
% In Fig.~\ref{fig:unknown}, the MC mean of the proposed estimate, computed over 100 runs, stabilizes after about $100$ iterations and converges toward MATLAB's value $\delta_\nu=0.3323$.
Using the same parameters as before ($N_{\mathrm{acc}}=10$, $N=9000$, $T_s=0.05$, and $\sigma_y^2=0.01$), the MC mean of the proposed estimate over 100 runs stabilizes after about $100$ iterations and converges in Fig.~\ref{fig:unknown} to MATLAB's value $\delta_\nu=0.3323$.
The dominant peak is attained at $\omega_{\mathrm{peak}}\approx 0.0126$ rad/sample, with $|p_0(\omega_{\mathrm{peak}})|\approx 0.644$, which places the maximizing frequency closer to nominal unity gain. From Vinnicombe's Riemann-sphere interpretation, this makes the geometric intuition behind the $\nu$-gap more directly relevant here; nevertheless, the larger value $\hat{\delta}_\nu=0.3323$ still indicates a non-negligible degree of closed-loop similarity. From a practical viewpoint, this suggests that a controller tuned for the GE 9F unit may also be applicable to the GE 6F unit. A formal guarantee, however, would require estimating the corresponding stability margin $b_{G_{9F},C}$, for example through a power-iteration-based estimate of $\|T(G_{9F},C)\|_\infty$. In that case, Proposition~\ref{prep:1} would apply whenever $b_{G_{9F},C}>0.3323$.

\subsection{Case Study: Poles in RHP}
We evaluate the $\nu$-gap metric on a discrete-time SISO case study with sampling time $T_s=1$\,s. Consider the stable plant
\begin{IEEEeqnarray}{rCl}
\label{eq:book}
G_1(z)=\frac{1-z^{-1}}{1-0.8z^{-1}},
\qquad
G_2(z)=1.8\,z^{-1}\;G_1(z).
\end{IEEEeqnarray}
Plant $G_2$ differs from $G_1$ by a gain of $1.8$ and an additional one-sample delay $z^{-1}$, which introduces extra phase lag. %and provides a compact test case where the $\nu$-gap estimation and the associated index (winding-number) condition can be examined jointly.
We can obtain explicit expressions for the index functions on the unit circle $z=e^{j\omega}$. 
Since the plants are SISO with real coefficients, $G_1(e^{-j\omega})=\overline{G_1(e^{j\omega})}$, and thus
\[
f_1(e^{j\omega})=1+G_1(e^{j\omega})G_1(e^{-j\omega})
=1+\big|G_1(e^{j\omega})\big|^2,
\]
\[
f_2(e^{j\omega})=1+G_2(e^{j\omega})G_1(e^{-j\omega})
=1+G_2(e^{j\omega})\,\overline{G_1(e^{j\omega})}.
\]
By~\eqref{eq:book}, we can simplify $f_2(e^{j\omega})$ as
\[
f_2(e^{j\omega})=1+1.8\,e^{-j\omega}\big|G_1(e^{j\omega})\big|^2.\vspace{-0.6em}
\]
Note that $f_1(e^{j\omega})$ is real and strictly positive for all $\omega$, hence $\mathrm{wno}(f_1)=0$. 
%By contrast, the one-sample delay introduces the factor $e^{-j\omega}$ in $f_2(e^{j\omega})$, whose rotation with $\omega$ may cause the curve to encircle the origin, yielding $\mathrm{wno}(f_2)\neq 0$.
In contrast, the factor $e^{-j\omega}$ induced by the one-sample delay rotates the complex term in $f_2(e^{j\omega})$ as $\omega$ varies, so the curve $f_2(e^{j\omega})$ may encircle the origin, yielding $\mathrm{wno}(f_2)\neq 0$.
\begin{figure}[h!]
\centering
    \vspace{-1.2 em}
    %\vspace{0.4 em}
\includegraphics[width=0.75\columnwidth]{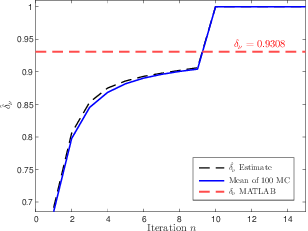}
\caption{Estimate of the $\nu$-gap $\hat{\delta}_\nu$ (black, dashed), average of 100 MC simulations (blue, solid), and actual $\delta_\nu$ (red, dashed).}
    \label{fig:disturbance} 
\end{figure}

Using $N_{\mathrm{acc}}=10$, $N=1000$, $T_s=1$, and $\sigma_y^2=0.01$, Fig.~\ref{fig:disturbance} shows that the power-iteration estimate initially approaches MATLAB's value $\delta_\nu=0.9308$. However, at $n=N_{\mathrm{acc}}$, Algorithm~\ref{alg:Ccheck} detects $\mathrm{wno}(f_1)\neq \mathrm{wno}(f_2)$, so the pair fails the set $\mathcal{C}$ condition. Hence, our procedure applies a hard stop and does not refine the estimate further.%, as the purpose of this example is to highlight the index-condition violation rather than continued numerical convergence.
%\subsection{Experimental Results}
%focusing on the accuracy of the estimated structured singular value and the impact of the underlying uncertainty structure $\Delta$ on its performance.

\section{Conclusion} \label{sec: conclusion}
We have presented an identification-free, data-driven method to estimate the Vinnicombe $\nu$-gap between two discrete-time SISO systems directly from time-domain input-output experiments. The proposed procedure combines iterative $\nu$-gap estimation with a data-driven admissibility check, providing a control-relevant measure of model mismatch for closed-loop analysis. Numerical studies, including an industrially motivated case, show close agreement with MATLAB’s \texttt{gapmetric} and correct detection of index-condition violations. Future work will address extension to MIMO systems.

\bibliography{References}

@book{aastrom_21feedback,
  title={Feedback Systems: An Introduction for Scientists and Engineers},
  author={{\AA}str{\"o}m, Karl Johan and Murray, Richard},
  year={2021},
  publisher={Princeton University Press}
}

@article{rowen_83,
    author = {Rowen, W. I.},
    title = {Simplified Mathematical Representations of Heavy-Duty Gas Turbines},
    journal = {Journal of Engineering for Power},
    volume = {105},
    pages = {865-869},
    year = {1983},
}

@inproceedings{khormali_15,
  title={Identification of an industrial gas turbine based on Rowen's model and using Multi-Objective Optimization method},
  author={Khormali, Aminollah and Yousefi, Iman and Yahyaei, Hassan and Aliyari, Sh Mahdi},
  booktitle={3rd RSI ICROM},
  pages={482--487},
  year={2015},
}

@article{tavakoli_09,
  title={An educational guide to extract the parameters of heavy duty gas turbines model in dynamic studies based on operational data},
  author={Tavakoli, Mohammad Reza Bank and Vahidi, Behrooz and Gawlik, Wolfgang},
  journal={IEEE Transactions on Power Systems},
  volume={24},
  number={3},
  pages={1366--1374},
  year={2009},
}

@inproceedings{NOX_96,
  title={Dry Low {NO}$_\text{x}$ Combustion Systems for {GE} Heavy-Duty Gas Turbines},
  author={Davis, L Berkley and Black, SH},
  year={1996},
  organization={ASME}
}

@inproceedings{venka_11,
  title={F-class {DLN} technology advancements: {DLN2}. 6+},
  author={Venkataraman, Krishna and Lewis, Skigh E and Natarajan, Jayaprakash and Thomas, Stephen R and Citeno, Joseph V},
  booktitle={Turbo Expo: Power for Land, Sea, and Air},
  volume={54631},
  pages={587--594},
  year={2011}
}

@article{welch_67,
  title={The use of fast {Fourier} transform for the estimation of power spectra: A method based on time averaging over short, modified periodograms},
  author={Welch, Peter},
  journal={IEEE Transactions on Audio and Electroacoustics},
  volume={15},
  number={2},
  pages={70--73},
  year={1967},
  publisher={IEEE}
}

@ARTICLE{guerrero_25,
  author={Guerrero, Margarita A. and Lakshminarayanan, Braghadeesh and Rojas, Cristian R.},
  journal={IEEE Control Systems Letters}, 
  title={Data-Driven Estimation of Structured Singular Values}, 
  year={2025},
  volume={9},
  number={},
  pages={1976-1981}}

@preamble{ " \newcommand{\noop}[1]{} " }

@inproceedings{vinnicombe_91,
  title={Structured uncertainty and the graph topology},
  author={Vinnicombe, Glenn},
  booktitle={Proceedings of the 30th IEEE CDC},
  pages={541--542},
  year={1991}
}

@inproceedings{vinnicombe_92,
  title={On the Frequency Response Interpretation of an Indexed {$\mathcal{L}_2$}-gap metric},
  author={Vinnicombe, Glenn},
  booktitle={American Control Conference},
  pages={1133--1137},
  year={1992}
}

@ARTICLE{sakkary_85,
  author={El-Sakkary, A.},
  journal={IEEE Transactions on Automatic Control}, 
  title={The gap metric: Robustness of stabilization of feedback systems}, 
  year={1985},
  volume={30},
  pages={240-247},}

@article{cantoni_17,
  title={Gap metric computation for time-varying linear systems on finite horizons},
  author={Cantoni, Michael and Pfifer, Harald},
  journal={IFAC-PapersOnLine},
  volume={50},
  number={1},
  pages={14513--14518},
  year={2017},
  publisher={Elsevier}
}

@book{desoer_09,
  title={Feedback Systems: Input-Output Properties},
  author={Desoer, Charles A and Vidyasagar, Mathukumalli},
  year={2009},
  publisher={SIAM}
}

@article{Megretski_97,
  title={System analysis via integral quadratic constraints},
  author={Megretski, Alexandre and Rantzer, Anders},
  journal={IEEE Transactions on Automatic Control},
  volume={42},
  pages={819--830},
  year={1997},
}

@article{wahlberg2010non,
  title={Non-parametric methods for $\mathcal{L}_2$-gain estimation using iterative experiments},
  author={Wahlberg, Bo and Syberg, M{\"a}rta Barenthin and Hjalmarsson, H{\aa}kan},
  journal={Automatica},
  volume={46},
  number={8},
  pages={1376--1381},
  year={2010},
  publisher={Elsevier}
}

@Book{Ljung:99,
  Title                    = {System Identification: Theory for the User, 2nd Edition},
  Author                   = {L. Ljung},
  Publisher                = {Prentice Hall},
  Year                     = {1999}
}

@article{rojas2012analyzing,
  title={Analyzing iterations in identification with application to nonparametric $\mathcal{H}_\infty$-norm estimation},
  author={C. R. Rojas and T. Oomen and H. Hjalmarsson and B. Wahlberg},
  journal={Automatica},
  volume={48},
  number={11},
  pages={2776--2790},
  year={2012}
}

@book{koch2022determining,
  title={Determining input-output properties of linear time-invariant systems from data},
  author={A. Koch},
  year={2022},
  publisher={Logos Verlag}
}

@article{koenings2017data,
  title={A data-driven computation method for the gap metric and the optimal stability margin},
  author={Koenings, T. and Krueger, M. and Luo, H. and Ding, S. X.},
  journal={IEEE Transactions on Automatic Control},
  volume={63},
  number={3},
  pages={805--810},
  year={2017},
  publisher={IEEE}
}

@inproceedings{muller2019gain,
  title={Gain estimation of linear dynamical systems using {Thompson sampling}},
  author={M. I M{\"u}ller and C. R. Rojas},
  booktitle={Proceedings of the 22nd International Conference on Artificial Intelligence and Statistics (AISTATS)},
  pages={1535--1543},
  year={2019}
}

@book{zhoudoyle,
author = {K. Zhou and J. C. Doyle and K. Glover},
title = {Robust and Optimal Control},
year = {1996},
publisher = {Prentice-Hall},
}

@book{golub2013matrix,
  title={Matrix Computations},
  author={G. H. Golub and C. F. {Van Loan}},
  year={2013},
  publisher={John Hopkins University Press}
}

%%%%%%%%%%%%%%%%%%%%%%%%%%%%%%%%%%%%%%%%%%%%%%%%%%%%%%%%%%%%%%%%%%%%%%%%%%%%%%%%

%%%%%%%%%%%%%%%%%%%%%%%%%%%%%%%%%%%%%%%%%%%%%%%%%%%%%%%%%%%%%%%%%%%%%%%%%%%%%%%%

\end{document}